\documentclass[11pt, reqno] {amsart}

\usepackage{amssymb}
\usepackage{hyperref} 
\usepackage{mathrsfs,bbold}
\usepackage{amsmath}
\usepackage{amsthm, mhequ}
\usepackage{enumerate}
\usepackage{dsfont}
\usepackage{color}
\usepackage[a4paper]{geometry}
\usepackage[utf8]{inputenc}
\usepackage{stmaryrd}
\usepackage{titlesec}
\usepackage{mathabx}  
\usepackage{appendix}
\usepackage{todonotes, fancyhdr}
\usepackage{chngcntr}
\usepackage{cancel}
\usepackage{breakurl}

\usepackage{setspace}
\renewcommand{\baselinestretch}{0.99}

	\definecolor{darkgreen}{rgb}{0.0, 0.5, 0.0}

\usepackage[all]{xy}
\numberwithin{subsection}{section}
\numberwithin{subsubsection}{subsection}
\numberwithin{equation}{section} 

\newenvironment{Dem}[1][\unskip]{%
    \begin{list}{\hspace{0.5cm}{\sf \textbf{Proof #1 --}}}{%
        \setlength{\topsep}{0pt}%
        \setlength{\leftmargin}{0pt}%
        \setlength{\rightmargin}{0pt}%
        \setlength{\listparindent}{0pt}%
        \setlength{\itemindent}{0pt}%
        \setlength{\parsep}{0pt}%
        \addtolength{\leftmargin}{20pt}%
        \addtolength{\rightmargin}{0pt}%
    } \item }{\hfill $\rhd$\end{list}\smallskip}

\renewcommand\thesection       {\arabic{section}}
\renewcommand\thesubsection    {\thesection{\boldmath $.$}\arabic{subsection}}
\renewcommand\thesubsubsection    {\thesection{\boldmath $.$}\arabic{subsection}{\boldmath $.$}\arabic{subsubsection}}

\titleformat{\section}[block]
{\filcenter\normalfont\sffamily\bfseries\Large}
{{\hspace{-0.7cm}}\thesection \hspace{0.2em} --\vspace{0.3cm}}{0.5em}{}

\titleformat{\subsection}[block]
{\filcenter\normalfont\sffamily\bfseries\large}  						  
{\hspace{-0.7cm}\thesubsection \hspace{0.5em} \vspace{0.3cm}}{.5em}{}  
\titlespacing{\subsection}{-0pc}{1.5ex plus .1ex minus .2ex}{0pc}

\titleformat{\subsubsection}[block]
{\filcenter\normalfont\sffamily\bfseries}					  
{\hspace{-0.7cm}\thesubsubsection \hspace{0.5em} \vspace{0.3cm}}{.5em}{}  
\titlespacing{\subsection}{-0pc}{1.5ex plus .1ex minus .2ex}{0pc}

\newtheoremstyle{mystyle}
{3pt}               
{3pt}               
{\it }                      
{}                      
{\sffamily\bfseries}             
{}                      
{0.5em}                 
{#1 #2{\Large$.$}  }

\theoremstyle{mystyle}

\newtheorem{thm}{Theorem}
\newtheorem*{thm*}{Theorem}

\newtheorem{cor}[thm]{\hspace{-0.15cm}  {Corollary} }
\newtheorem{prop}[thm]{\hspace{-0.13cm} {Proposition}}
\newtheorem{rem}[thm]{\hspace{-0.15cm} {Remark}}

\newtheoremstyle{mystyle2}
{3pt}               
{3pt}               
{\it }                      
{}                      
{\sffamily\bfseries}             
{}                      
{0.5em}                 
{\llap{#2 }#1{\hspace{0.2cm}--}}

\theoremstyle{mystyle2}

\newtheorem*{definition*}{Definition}
\newtheorem*{theorem*}{Theorem}
\newtheorem*{Remark*}{Remark}
\newtheorem*{lem*} {Lemma}
\newtheorem*{defn*} {Definition}
\newtheorem*{prop*} {Proposition}
\newtheorem*{cor*} {Corollary}

\renewcommand{\epsilon}{\varepsilon}

\newcommand{\mcB}{\mathcal{B}} 
 
\newcommand{\mcD}{\mathcal{D}}

\def\Labe{\mathfrak{e}}

\def\Labn{\mathfrak{n}}

\def\Labhom{\mathfrak{t}}

\def\Lab{\mathfrak{T}}
\def\CT{\mathcal{T}}

\def\CI{\mathcal{I}}
\def\N{\mathbb{N}}
\def\one{\mathbf{1}}
\def\s{\frak{s}}
\def\root{\mathrm{root}}
\def\nonroot{\mathrm{non-root}}

\def\id{\mathrm{id}}

\newcommand{\RR}{\mathbb{R}}

\DeclareMathAlphabet{\mcb}{U}{BOONDOX-calo}{m}{n}
\SetMathAlphabet{\mcb}{bold}{U}{BOONDOX-calo}{b}{n}
\DeclareFontFamily{U}{mathx}{\hyphenchar\font45}
\DeclareFontShape{U}{mathx}{m}{n}{
      <5> <6> <7> <8> <9> <10>
      <10.95> <12> <14.4> <17.28> <20.74> <24.88>
      mathx10
      }{}
\DeclareSymbolFont{mathx}{U}{mathx}{m}{n}
\DeclareMathSymbol{\bigtimes}{1}{mathx}{"91}
\def\bigtimes{\times}

\begin{document}

\begin{center}
{\LARGE\sffamily{Renormalised singular stochastic PDEs   \vspace{0.5cm}}}
\end{center}

\begin{center}
{\sf I. BAILLEUL} \& {\sf Y. BRUNED}
\end{center}

\vspace{1cm}

\begin{center}
\begin{minipage}{0.8\textwidth}
\renewcommand\baselinestretch{0.7} \scriptsize \textbf{\textsf{\noindent Abstract.}} Extended decorations on naturally decorated trees were introduced in the work of Bruned, Hairer and Zambotti on algebraic renormalization of regularity structures to provide a convenient framework for the renormalization of systems of singular stochastic PDEs within that setting. This non-dynamical feature of the trees complicated the analysis of the dynamical counterpart of the renormalization process. We provide a new proof of the renormalised system by-passing the use of extended decorations and working for a large class of renormalization maps, with the BPHZ renormalization as a special case. The proof reveals important algebraic properties connected to preparation maps.
\end{minipage}
\end{center}

\vspace{0.6cm}

\section{Introduction}
\label{SectionIntro}

We consider systems of parabolic equations involving possibly different second order elliptic linear differential operators with constant coefficients $L_1,\dots, L_{k_0}$
\begin{equation} \label{main_equation}
(\partial_t - L_i)u_i = F_i(u, \nabla u)\xi,   \qquad (i=1\dots k_0),
\end{equation}
with $u:=(u_1,\dots,u_{k_0})$ and each $u_i$ taking values in a finite dimensional space $\RR^{d_i}$, with $F_i(u, \nabla u)\in L(\RR^{n_0},\RR^{d_i})$, for each $u$, and $\xi = (\xi_1,\dots,\xi_{n_0})$ an $n_0$-dimensional spacetime `noise'. The unknown is defined on positive times, with a given time $0$ initial value. Some of the components of $\xi$ may be regular, or even constant, functions like in the generalized (KPZ) equation 
$$
(\partial_t-\partial_x^2) u = f(u)\xi_1 + g(u)\vert\partial_x u\vert^2,
$$
where $u$ is a real-valued function, $\xi=(\xi_1,\xi_2)$, with $\xi_1$ a one-dimensional spacetime white noise and $\xi_2=1$. The foundations of regularity structures theory are contained in M. Hairer's groundbreaking work \cite{reg} and his subsequent works \cite{BHZ, CH16, BCCH} with Bruned, Chandra, Chevyrev and Zambotti. We refer the reader to Friz and Hairer's book \cite{FrizHai} for a gentle introduction to the subject, to Bruned, Hairer and Zambotti's short review \cite{EMS}, and to Bailleul and Hoshino's work \cite{BaiHos} for a complete concise account of the analytic and algebraic sides of the subject. Possible solutions to a given equation/system are defined by their local behaviour
$$
u(\cdot) \simeq \sum_\tau u_\tau(x)({\sf \Pi}_x\tau)(\cdot),
$$
in terms of reference functions/distributions $({\sf \Pi}_x\tau)(\cdot)$, indexed by all state space points $x$ and a finite collection of symbols $\{\tau\}$. The reconstruction theorem ensures that one defines indeed in this way a unique function/distribution when the coefficient $\big\{u_\tau(x)\big\}$  form a consistent family, encoded in the notion of modelled distribution. The minimum set of symbols $\{\tau\}$ needed to give local expansions of possible solutions to singular PDEs have a natural combinatorial structure that comes with the naive Picard formulation of the equation and the very fact that they are used to build local expansion devices. Regularity structures are the appropriate abstraction of these combinatorial structures and models on regularity structures the appropriate analogue of local expansion devices.

Under proper subcriticality conditions on a given system of  singular stochastic PDEs, one can build a regularity structure $\mathscr{T}$ for it, and given a model ${\sf M=(g ,\Pi)}$ on $\mathscr{T}$, one can recast the system \eqref{main_equation} as a fixed point problem for a modelled distribution ${\sf u}\in \mcD^\gamma(T,{\sf g})$, of regularity $\gamma$, for a model-dependent equation in the space $\mcD^\gamma(T,{\sf g})$, for a well-chosen positive regularity exponent $\gamma$. Such functions take values in a finite dimensional linear subspace $T_{<\gamma}$ of the linear space $T$. Obtaining a fixed point usually requires adjusting a parameter to get a contractive map on a proper functional space. This typically gives well-posedness results in small time or small parameter. A proper distribution/function on the state space where $\sf u$ is defined is obtained by applying to $\sf u$ the reconstruction operator $\sf R^M$ associated with the model $\sf M$. While the modelled distribution $\sf u$ solves a dynamical equation, counterpart of system \eqref{main_equation}, the possibility to give a dynamical description of its reconstruction $\sf R^Mu$, and to relate it to the formal equation \eqref{main_equation}, depends on the model. Denote by $\zeta=(\zeta_1,\dots,\zeta_{n_0})$ the symbol in $\mathscr{T}$ used to represent the noise $\xi$.

The use of admissible models ensures that
\begin{equation} \label{EqReconstructionDynamics}
u_i = K_i * {\sf R^M}\big({\sf F}_i({\sf u})\,\zeta \big) + K_i\big(u_i(0)\big),   \qquad (i=1\dots k),
\end{equation} 
where $K_i$ is the heat kernel of the operator $L_i$, the value at time $0$ of $u_i$ is $u_i(0)$, and the symbol $*$ stands for spacetime convolution. The functions ${\sf F}_i$ are the natural extensions of the functions $F_i$ to the space of modelled distributions. If the noise $\xi$ is smooth and one uses the canonical admissible model $\Theta$ sending the noise symbol $\zeta$ to $\xi$, its reconstruction operator happens to be multiplicative and system \eqref{EqReconstructionDynamics} turns out to be equivalent to system \eqref{main_equation}. The canonical admissible model is no longer well-defined if the noise is not sufficiently regular, which is the case of interest. In that case, one can only build random models, using probability tools, when the noise itself is random and satisfies some mild conditions detailed in Chandra and Hairer's work \cite{CH16}. These admissible models $\sf M$ are limits in probability of admissible models ${\sf M}^\epsilon = ({\sf g}^\epsilon, {\sf \Pi}^\epsilon)$ for which 
$$
{\sf \Pi}^\epsilon = {\sf \Theta}^\epsilon\circ R^\epsilon,
$$
where ${\sf \Theta}^\epsilon$ is the naive interpretation operator on symbols mapping the noise symbol $\zeta$ to a regularized version $\xi^\epsilon$ of $\xi$ that respects the parabolic scaling of the equation, and $R^\epsilon$ is a deterministic linear map on the finite dimensional linear space $T_{<\gamma}$, diverging as the regularization parameter $\epsilon$ goes to $0$. It is then no longer clear at all that the ${\sf M}^\epsilon$-reconstruction $u^\epsilon$ of the solution ${\sf u}^\epsilon$ to the above mentioned fixed point problem in $\mathcal{D}^\gamma(T,{\sf g}^\epsilon)$ is the solution of a PDE involving the regularized noise $\xi^\epsilon$. As the model ${\sf R}^{{\sf M}^\epsilon}$ takes values in the space of continuous functions, its reconstruction operator ${\sf R}^{{\sf M}^\epsilon}$ satisfies
$$
\big({\sf R}^{{\sf M}^\epsilon} {\sf v}\big)(x) = \big({\sf \Pi}^\epsilon_x{\sf v}(x)\big)(x),
$$
for all modelled distributions $\sf v$ of positive regularity. The possibility to turn the non-autonomous dynamics \eqref{EqReconstructionDynamics} for $u^\epsilon$ and ${\sf u}^\epsilon$ into an autonomous dynamics for $u^\epsilon$ depends then on our ability to compute effectively the recentered renormalized interpretation operator ${\sf \Pi}^\epsilon_x$ associated with $({\sf g}^\epsilon, {\sf \Pi}^\epsilon)$. This is a non-elementary matter. Hairer used in his seminal work \cite{reg} the fact that one has for the $\Phi^4_3$ and $2$-dimensional generalized (PAM) equations
\begin{equation} \label{EqCondition1}
\big({\sf \Pi}^\epsilon_x\tau\big)(x) = \big({\sf \Theta}^\epsilon_x(R^\epsilon\tau)\big)(x), \qquad \forall\,\tau\in T,\;\forall\,x,
\end{equation}
for the natural choice of decorated trees $\tau$ associated with these equations, to deal by hand with these equations. Such a property of sets of natural trees associated with singular PDEs was later proved to hold as well for the $\sin$-Gordon equation \cite{CHS}, the generalized (KPZ) equations and the $\Phi^4_{4-\delta}$ equation \cite{BGHZ, BCCH}, and the $2$-dimensional Yang-Mills equation \cite{CCHS}. The stronger property
\begin{equation} \label{EqCondition2}
{\sf \Pi}^\epsilon_x\tau= \Theta^\epsilon_x(R^\epsilon\tau), \qquad \forall\,\tau\in T,\;\forall\,x,
\end{equation}
holds for the list of trees that comes from the Picard development of solutions of the (KPZ) and generalized (PAM) equations. However, not all subcritical singular PDEs satisfy either of these properties. The introduction in Bruned, Hairer and Zambotti's work \cite{BHZ} of extended decorations on the set of decorated trees was motivated by the desire to set a framework where this identity holds true for a whole class of equations. 

Bruned, Chandra, Chevyrev and Hairer showed in \cite{BCCH} that one can run within this framework a clean analysis of the dynamics of $u^\epsilon$, and that $u^\epsilon=\big(u^\epsilon_1,\dots,u^\epsilon_k\big)$ is solution of the system

\begin{equation} \label{EqRenormalisedEquation}
(\partial_t - L_i)u_i^\epsilon = F_i(u^\epsilon,\nabla u^\epsilon)\xi^\epsilon + \sum_{\tau\in \mathcal{B}^-\backslash\{\textbf{\textsf{1}}\}} \ell^\epsilon(\tau)\,\frac{F_i(\tau)(u^\epsilon, \nabla u^\epsilon)}{S(\tau)},   \qquad (i=1\dots k_0),
\end{equation}
with additional explicit counterterms $F_i(u^\epsilon, \nabla u^\epsilon)$ depending on $u^\epsilon$ and its derivative, possibly, and where $\ell^\epsilon(\tau)$ are renormalization constants indexed by a finite set of decorated trees $\mathcal{B}^-$ containing an element \textbf{\textsf{1}}, diverging as $\epsilon$ goes to $0$, and $S(\tau)$ is a symmetry factor. The terms $F_i(\tau)(\cdot)/S(\tau)$ are the counterparts of the coefficients used to describe $B$-series in numerical analysis. (Decorated
trees and the same type of coefficients have been used recently to describe a numerical scheme
at low regularity for dispersive equations in \cite{BS}.) This comparison makes less surprising the crucial role plaid by pre-Lie structures in the analysis of singular PDEs and the fact that the preceding terms satisfy some crucial morphism property for a (multi-)pre-Lie structure that was first introduced by Bruned, Chandra, Chevyrev and Hairer in \cite{BCCH}. (That such structures have a role to play in these questions was first noticed in a rough paths setting in Bruned, Chevyrev, Friz and Preiss' works \cite{BCF, BCFP}.) Equation \eqref{EqRenormalisedEquation} is called the {\it renormalized equation}. At the algebraic level, identity \eqref{EqCondition2} reflects a co-interaction between two Hopf algebras, whose use in \cite{BCCH} for deriving the renormalized equation for a large class of singular PDEs was instrumental  -- see e.g. Section 5 of \cite{BaiHos} for the core points of this co-interaction.

\smallskip

We show in this work that none of the idendities \eqref{EqCondition1} and \eqref{EqCondition2} is actually needed to get back the renormalized equation and that one can run the analysis in the natural space of trees with no extended decorations.

\smallskip

A systematic renormalization procedure was designed by Bruned, Hairer and Zambotti in \cite{BHZ} and proved to provide converging renormalized models by Chandra and Hairer in \cite{CH16}. It is named `BPHZ renormalization', after similar renormalization procedures introduced by Bogoliubov and Parasiuk for the needs of quantum field theory in the mid 50's, improved among others by Hepp and Zimmermann. Its regularity structures counterpart is subtle as renormalization needs to cope well with the recentering properties of the model. This compatibility between recentering and renormalization structures is encoded at an algebraic level in the above mentioned co-interaction of two Hopf algebras. (Such a co-interaction has been observed by Chartier, E. Hairer, Vilmart and Calaque, Ebrahimi-Fard, Manchon's works \cite{CHV10,CEM} in the simpler context of the Butcher-Connes-Kreimer and the extraction-contraction Hopf algebras.) We consider here a larger class of renormalization procedures introduced in \cite{BR18}, containing the BPHZ renormalization as an element. They are built from special linear maps $R : T\rightarrow T$, to which one can associate a family of multiplicative operators $\Pi_x^{M^\circ}$ from $T$ to the space of smooth functions, and a smooth admissible model ${\sf M=(g,\Pi)}$ on $T$ such that its associated reconstruction operator $\sf R^M$ on modelled distributions of positive regularity factorizes through a multiplicative map
$$
\big({\sf R^Mv}\big)(x) = \Big(\Pi^{M^\circ}{\sf v}(x)\Big)(x).
$$
This brings back the core problem to understanding the action of $R$ on the lift to the regularity structure of the equation. The key point is then a right morphism property of $R$ for an $R$-independent product $\star$ introduced by Bruned and Manchon in \cite{deformation} as the dual of the deformed Butcher-Connes-Kreimer coproduct used in \cite{BHZ}.

\smallskip

The remainder of this work is organised as follows. In Section~\ref{SectionDecoDeformedPreLie}, we recall basics on decorated trees and deformed pre-Lie structures using mainly the formalism developed in Bruned and Manchon's work \cite{deformation}. The main new result in this section is Proposition \ref{star_morphism}, which establishes a morphism property for $ F $ with respect to the $\star$-product. We introduce good multi-pre-Lie morphisms and (strong) preparation maps $R$ in Section \ref{SectionLocalProduct}, and show in Proposition \ref{prop_equivalence} that adjoints $R^*$ of strong preparation maps are right-morphisms for the $\star$-product. Proposition \ref{PropCombined} combines this result with Proposition \ref{star_morphism} to provide a clear understanding of the action of $R^*$ on $F$. The main result of the section is Theorem \ref{adjoint_M}. It states that strong preparation maps $M$ are good multi-pre-Lie morphisms. Corollary \ref{R_translation} states that the set of good multi-pre-Lie morphisms is in bijection with the set of preparation maps in a rough paths setting; a remarkable result. We prove the main result of the paper, Theorem \ref{ThmMain}, in Section \ref{SectionShortProof}. It shows that the reconstruction of the solution to the singular PDE in the space of modelled distributions associated with our class of renormalization maps is actually the solution of an explicit PDE of the form
$$
(\partial_t-L_i)u_i = F_i(u,\nabla u)\xi + \sum_{1\leq l\leq n_0}F_i\big((R^*-\textrm{Id})\zeta_l\big)(u,\nabla u)\xi_l,\qquad (1\leq i\leq k_0).
$$
One gets back a system of the form \eqref{EqRenormalisedEquation} for maps $R^*$ fixing symbols $\zeta_l$ corresponding to non-constant noises.

\bigskip

\noindent \textbf{\textsf{Notation --}} The letter $\zeta$ will be used exclusively for the noise symbol in a regularity structure. The letters $\sigma, \tau, \mu, \nu$ will denote (decorated) trees.

\section{Decorated trees and pre-Lie products}
\label{SectionDecoDeformedPreLie}

Recall system \eqref{main_equation} with its noises $\xi_1,\dots,\xi_{n_0}$ and its operators $L_1,\dots, L_{k_0}$. Let
$$
\mathfrak{T}^-=\big(\Labhom^-_1,\dots,\Labhom^-_{n_0}\big), \;\textrm{and}\; \mathfrak{T}^+=\big(\Labhom_1^+,\dots,\Labhom_{k_0}^+\big)
$$ 
be finite sets representing noise types and operator types, respectively. Denote by 
$$
\mathfrak{T} := \mathfrak{T}^-\cup\mathfrak{T}^+
$$
the set of all types. We consider decorated trees $(\tau,\Labn,\Labe)$ where $\tau$ is a non-planar rooted tree with node set $N_\tau$ and edge set $E_\tau$. The maps $ \Labn : N_\tau \rightarrow \N^{d+1} $, and $\Labe=\big(\frak{t}(\cdot),\frak{p}(\cdot)\big): E_\tau \rightarrow \Lab \times \N^{d+1}$, are node decorations and edge decorations, respectively. The $\N^{d+1}$-part $\frak{p}(e)$ in the edge decoration of an edge $e$ encodes possible derivatives acting on the operator associated with the given edge type $\frak{t}(e)$. We will frequently abuse notations and simply denote by $\tau$ a decorated tree, using a symbolic notation. 
\begin{enumerate}
   \item[$\bullet$] An edge decorated by  $  (\Labhom,p) \in \Lab \times \N^{d+1} $  is denoted by $ \CI_{(\Labhom,p)} $. The symbol $  \CI_{(\Labhom,p)} $ is also viewed as  the operation that grafts a tree onto a new root via a new edge with edge decoration $ (\Labhom,p) $
    \item[$\bullet$] A factor $ X^k $ encodes a single node  $ \bullet^{k} $ decorated by $ k \in \N^{d+1} $. Denote by $\{e_1, \ldots, e_{d+1}\}$ form the canonical basis of $ \N^{d+1} $. For $1\leq i\leq d+1$, write $ X_i $ for $ X^{e_i} $. The element $ X^0 $ is identified with $ \one $.
\end{enumerate}
We require that every decorated tree $ \tau $ contains at most one edge decorated by $ (\Labhom,p) $ with $ \Labhom \in \Lab^- $ and any $p\in\N^{d+1}$, at each node. This encodes the fact that no product of two noises are involved in the analysis of the system \eqref{main_equation}. We suppose that these edges lead directly to leaves; we denote them by $\zeta_l$, for $1\leq l\leq n_0$; by convention $ \zeta_0 $ is equal to $ \one $.  Any decorated tree $ \tau $ has a unique decomposition
 \[
	\tau =  X^{k} \zeta_l \prod_{i=1}^{n} \CI_{a_i}(\tau_i) ,
 \]
 where $\prod_i$ is the tree product, the $\tau_i$ are decorated trees and the $a_i$ belong to $ \Lab^+ \times \N^{d+1}$, so no factor in the product is a noise symbol $\zeta_{l'}$. The algebraic symmetry factor $S(\tau)$ of a decorated tree $\tau= X^{k} \zeta_l \ \prod_{j=1}^m \mathcal{I}_{a_j}(\tau_j)^{\beta_j}$ is defined grouping terms uniquely in such a way that $(a_i,\tau_i) \neq (a_j,\tau_j)$ for $i \neq j$, and setting inductively
\begin{equation*}
S(\tau) = k!\,
\bigg(
\prod_{j=1}^{m}
S(\tau_{j})^{\beta_{j}}
\beta_{j}!
\bigg)\;.
\end{equation*} 
A {\it planted tree} is a tree of the form $\CI_a(\sigma)$, for a decorated tree $\sigma$ and $a\in\mathfrak{T}^+\times\mathbb{N}^{d+1}$; we denote by $\CI(T)$ the set of planted trees. We define an inner product on the set of all decorated trees setting for all $\sigma, \tau \in\CT$
\begin{align*} \label{inner_product}
\langle \sigma ,  \tau \rangle := S(\tau)\,\textbf{\textsf{1}}_{\sigma=\tau}.
\end{align*}
We also set 
$$
\langle \sigma_1\otimes\sigma_2 ,  \tau_1\otimes\tau_2 \rangle := \langle \sigma_1 , \tau_1 \rangle\,\langle \sigma_2 , \tau_2 \rangle.
$$
The linear span of decorated trees will be denoted by $T$. Note here that such trees are also useful to describe numerical schemes for dispersive equations with low regularity initial condition \cite{BS}.

\smallskip

We now associate numbers to decorated trees. Fix a scaling $ \frak{s} \in \N^{d+1} $ and a map 
$$
|\cdot|_{\s} : \Lab \rightarrow \mathbb{R},
$$ 
which is negative on the noise types $ \Lab_-$ and positive on the operator types $\Lab_+$. This map accounts for the regularity of the noises and the gain of regularity of the heat kernels $K_i$, encoded in Schauder-type estimates they satisfy. We extend the map $|\cdot|_{\s}$ to $\mathfrak{T}\times\N^{d+1} $ setting 
$$
 |p|_{\s}  := \sum_{i=1}^{d+1} \s_i p_i, \qquad \textrm{and}\qquad |(\frak{t},p)|_{\s} := |\frak{t}|_{\s} +  |p|_{\s}, \qquad \textrm{for}\;k\in\N^{d+1}.
 $$ 
 The degree of a decorated rooted tree $(\tau, \Labn,\Labe)$  is defined by
\begin{equs}
\deg(\tau, \Labn,\Labe) := \sum_{v \in N_{T}} \big|\Labn(v)\big| _{\s}+ \sum_{e \in E_{T}} \big|\frak{t}(e)\big| _{\s} - \big|\frak{p}(e)\big| _{\s}.
\end{equs}
(`Degree' is called `homogeneity' in Hairer's work \cite{reg}.) We use the degree to introduce the space of `positive' decorated trees $T_+$. It is the linear span of trees of the form $X^k \prod_{i=1}^{n} \CI_{a_i}(\tau_i)$, where $\deg( \CI_{a_i}(\tau_i)) > 0$ and $k\in\N^{d+1}$. We also consider the linear space $ T_- $ spanned by the decorated trees with negative degree, and denote by $\mathbb{R}[T_-] $ the linear space spanned by forests of trees in $T_-$.

\smallskip

Given $k\in\N^{d+1}$ denote by $\uparrow_v^k$ the derivation on decorated trees that adds $k$ to the decoration at the node $v$. We introduce a family of pre-Lie products of grafting type setting for all decorated trees $\sigma, \tau\in T$, and $a\in\frak{L}\times\N^{d+1}$,
\begin{equation*}
\sigma \curvearrowright_a \tau := \sum_{v\in N_{\tau}}\sum_{m\in\N^{d+1}}{\Labn_v \choose m} \,\sigma  \curvearrowright^v_{a-m}(\uparrow_v^{-m} \tau),
\end{equation*}
where $ \Labn_v $ is the decoration at the node $ v $, and $\curvearrowright^v_{a-m}$ grafts $ \sigma $ onto $\tau$ at the node $ v $ with an edge decorated by $a -m$. (For $a=(\frak{t},p)$, one writes $a-m$ for $(\frak{t},p-m)$.) This formula requires that $\sigma=\textbf{\textsf{1}}$ is the only argument accepted on the left of the grafting operation when $a\in\frak{T}^-\times\N^{d+1}$, since edges of noise type have no other arguments. The above sum is finite due to the binomial coefficient $ {\Labn_v \choose m} $, which is equal to zero if $m$ is greater than $ \Labn_v $, by convention. The pre-Lie products $\curvearrowright_a$ are non-commutative; they were first introduced in Bruned, Chandra, Chevyrev and Hairer's work \cite{BCCH}. We recall one universal result that we will use in the sequel; it was first established in Corollary 4.23 of \cite{BCCH}. It can be viewed as an extension of the universal result of Chapoton-Livernet \cite{ChaLiv} on pre-Lie algebras. (Such a result becomes immediate when one constructs $\curvearrowright_a$  as a deformation, as in Section 2.1 of \cite{deformation}. See also Foissy's work \cite{F2018} for the case with no deformation.)

\smallskip

\begin{prop} \label{freeness} 
The space $T$ is freely generated by the elements $\Big\{ X^k \zeta_l; \, 1\leq l\leq n_0, \ k \in \N^{d+1}\Big\}$ and the operations $\big\{\hspace{-0.08cm}\curvearrowright_a\,;\,1\leq a\leq k_0\big\}$.
\end{prop}

\smallskip

We define a product 
$$
\curvearrowright\,: \CI(T)\times T\rightarrow T
$$ 
setting for all $a\in\mathfrak{T}\times\N^{d+1}, \sigma, \tau\in T$, with the appropriate restriction on $\sigma$ if $a\in\mathfrak{T}^-\times\N^{d+1}$,
\begin{equs} \label{def_pre_Lie}
\CI_a(\sigma ) \, \curvearrowright \, \tau := \sigma \, \curvearrowright_a \, \tau,
\end{equs}
We extend this product to a product of planted trees, $\prod_{i=1}^n \CI_{a_i}(\sigma_i ) \,\curvearrowright \, \tau$, by grafting each tree $ \sigma_i $ on $ \tau $ along the grafting operator $\curvearrowright_{a_i} $, independently of the others -- we allow here one of the $a_i$ to be an element of $\frak{T}^-\times\N^{d+1}$, so the product contains in that case (only) one noise.(The trees $\sigma_i$ are only grafted on $\tau$, not on one another.) Following Bruned and Manchon's construction in \cite{deformation}, for  $ B \subset N_{\tau} $, consider the derivation map $ \uparrow^{k}_{B}$ defined as
\begin{equs}
\uparrow^{k}_{B} \tau = \sum_{\sum_{v \in B} k_v  = k } \prod_{v \in B}\uparrow_v^{k_v}  \tau.
\end{equs}
and set, as a shorthand notation for later use,
$$
\uparrow^k \tau := \uparrow^k_{N_\tau} \tau.
$$
We define the product
$$
\star : T\times T\rightarrow T,
$$
for all $\sigma = X^k \prod_{i} \CI_{a_i}(\sigma_i) \in T$ and $\tau \in T$, by the formula
\begin{equs}
\sigma \star \tau :=  \,\uparrow^k_{N_\tau} \left( \prod_{i} \mathcal{I}_{a_i}(\sigma_i)\curvearrowright \tau  \right)
\end{equs}
One has for instance
$$
X^k\,\zeta_l\prod_{i=1}^n\CI_{a_i}(\tau_i) = \left(X^k\prod_{i=1}^n\CI_{a_i}(\tau_i)\right)\star \zeta_l.
$$
It has been proved in Section 3.3 of \cite{deformation} that this product is associative; this can be obtained by applying the Guin-Oudom procedure \cite{Guin1,Guin2} to a well-chosen pre-Lie product. When $ \sigma \in T_+ $ and $\tau, \mu\in T$, one has from Theorem 4.2 in \cite{deformation}
\begin{equs} \label{def_coaction}
\langle  \sigma \star \tau ,  \mu \rangle  = \langle  \tau \otimes \sigma ,  \Delta \mu \rangle 
\end{equs}
where 
$$
\Delta : T \rightarrow T \otimes T_+ 
$$ 
is a co-action first introduced in Hairer' seminal work \cite{reg} -- see also \cite{BHZ} and \cite{BaiHos}, where it plays a prominent role. So the restriction to $T^+\times T$ of the product $\star$ is the $\langle\cdot,\cdot\rangle$-dual of the splitting map $\Delta$. (Note here that our product $\star$ corresponds to the $\star_2$ product in \cite{deformation}.)

\smallskip

With a view on the system \eqref{main_equation} of singular PDEs, assume we are given a family 
$$
(F_k^l)_{1\leq k\leq k_0, 1\leq l\leq n_0}
$$ 
of functions of abstract variables $Z_a$ indexed by $a\in\frak{L}^+\times\N^{d+1}$. These variables account for the fact that the nonlinearities in \eqref{main_equation} may depend on $u$ and its derivatives -- only $u$ and $\nabla u$ for \eqref{main_equation}, but we could also consider systems where the differential operators $L_i$ have order higher than $2$; in which case the nonlinearities could depend on $u$ and all its $k$-th derivatives, for $k$ up to the order of $L_i$ minus $1$. The different components of $u$ are indexed by $1\leq i\leq k_0$, and its derivatives by $\N^{d+1}$ -- with $d$ space dimensions and one time dimension. We define in the usual way partial derivatives $D_a$ in the variable $Z_a $, and set for all $k\in\N^{d+1}$
$$ 
\partial^k := \sum_a Z_{a+k}\,D_a.
$$
We define inductively a family $F = (F_i)_{1\leq i\leq k_0}$ of  functions of the variables $Z_a$, indexed by $T$, setting for $ \tau = X^{k} \zeta_{l}\prod_{j=1}^{n} \CI_{a_j}(\tau_j) $, with $a_j = (\Labhom_{l_j},k_j)$, for all $1\leq i\leq k_0$,

\begin{equation} \label{def_upsilon}
\begin{split}
F_i(\zeta_l) := F^{l}_i, \qquad F_i(\tau) := \partial^k D_{a_1} ... D_{a_n} F_i (\zeta_l)\,\prod_{j=1}^n F_{l_j}(\tau_j).
\end{split} \end{equation}
The next statement is a morphism property of $F$ for the product $\star$, as a function on $T$.

\medskip

\begin{prop} \label{star_morphism}
For every $1\leq i\leq k_0$, for every $ \sigma\in T$, and $\tau = X^{k} \prod_{j=1}^{n} \CI_{a_j}(\tau_j) \in T$ with $ a_j = (\Labhom_{l_j}, k_j) $, one has
\begin{equs}\label{star2_morphism}
F_i\bigg(\Big\{X^{k} \prod_{j=1}^{n} \CI_{a_j}(\tau_j)\Big\} \star \sigma\bigg) = \partial^k D_{a_1} ... D_{a_n} F_i(\sigma)\,\prod_{j=1}^n F_{l_j}(\tau_j).
\end{equs}
\end{prop}

\smallskip

\begin{Dem}
We proceed by induction on $ \sigma $. The case $ \sigma = \zeta_l $ is part of the definition \eqref{def_upsilon}. For $ \sigma = X^{m} \zeta_l $, we have
\begin{equs}
X^{k} \prod_{j=1}^{n} \CI_{a_j}(\tau_j) \star \sigma = X^{k+m -  \sum_j \ell_j} \sum_{\ell \in (\N^{d+1})^{n}} {m \choose \ell} \prod_{j=1}^n \CI_{a_j-\ell_i}(\tau_j).
\end{equs}
Using the fact that
\begin{equs} \label{commutation derivatives}
\sum_{\ell=(\ell_1,\dots,\ell_n) \in (\N^{d+1})^{n}} {m \choose \ell} \partial^{m-\ell_j} D_{a_j-\ell_j} = D_{a_j} \partial^{m},
\end{equs}
one has
\begin{equs}
F_{\Labhom}\bigg(X^{k} \prod_{j=1}^{n} \CI_{a_j}(\tau_j) \star \sigma\bigg) & = \prod_{j=1}^{n}F_{l_j}(\tau_j) \sum_{\ell \in (\N^{d+1})^{n}} {m \choose \ell} \partial^{k+m- \sum_j \ell_j}   \prod_{j'} D_{a_{j'}-\ell_{j'}} F_i (\zeta_l)   \\
& = \prod_{j=1}^{n}F_{l_j}(\tau_j) \,\partial^{k}   \prod_{j'} D_{a_{j'}} D^{m}F_i(\zeta_l)   \\ 
&= \prod_{j=1}^{n}F_{l_j}(\tau_j) \,\partial^{k}   \prod_{j'} D_{a_{j'}} F_i(X^{m}\zeta_l) 
\end{equs}
Then, we assume that $ \sigma  = \CI_{b_1}(\sigma_1)\curvearrowright\sigma_2  = \CI_{b_1}(\sigma_1) \star \sigma_2 $ with $ b_1 = (\Labhom_b,k_b) $. One has from the associativity of $\star$
\begin{equs}
& \tau  \star \sigma  = \big( \tau \star \CI_{b_1}(\sigma_1) \big) \star \sigma_2
\\ & = \left\{ \sum_{I \subset \lbrace 1,...,n \rbrace}\sum_{k_1 + k_2 = k}  X^{k_1} \prod_{j \in I} \CI_{a_j}(\tau_j)  \CI_{b_1}\left(  X^{k_2} \prod_{j'\in \lbrace 1,...,n\rbrace \setminus I}  \CI_{a_{j'}}(\tau_{j'}) \star \sigma_1 \right) \right\} \star   \sigma_2
\end{equs}
We apply the induction hypothesis to get
\begin{equs}
F_i(\tau \star \sigma) & = \sum_{I \subset \lbrace 1,...,n \rbrace}\sum_{k_1 + k_2 = k}  \prod_{j =1}^{n} F_{l_j}(\tau_j)   \partial^{k_2} \prod_{j\in \lbrace 1,...,n\rbrace \setminus I} D_{a_j} F_b(\sigma_1) \,\partial^{k_1}\prod_{j'} D_{a_{j'}} D_{b_1} F_i(\sigma_2)
\end{equs}
Using the fact that $ \partial^{k} $ and $ D_{a_i} $ satisfy the Leibniz rule one then gets
\begin{equs}
F_i(\tau \star \sigma) & =  \prod_{j=1}^n F_{l_j}(\tau_j) \, \partial^{k} \prod_{j' =1}^n D_{a_{j'}} F_b(\sigma_1) D_{b_1}  F_i(\sigma_2)   \\
 &= \prod_{j =1}^n F_{l_j}(\tau_j) \, \partial^{k} \prod_{j' =1}^n D_{a_{j'}} F_i\Big(\CI_{b_1}(\sigma_1) \curvearrowright\sigma_2\Big)   \\ 
 &= \prod_{j =1}^n F_{l_j}(\tau_j)  \, \partial^{k} \prod_{j' =1}^n D_{a_{j'}} F_i(\sigma).
\end{equs}
which allows us to conclude the proof.
\end{Dem}

\medskip

Identity \eqref{commutation derivatives} was first noticed in the proof of Proposition 30 in (the first version of) Bailleul and Hoshino's work \cite{BaiHos}. It lead the authors to a simple proof of the fact that for $a=(\frak{t}_j,p_a)$ and all $i$
\begin{equs}
F_i\big(\mathcal{I}_{a}( \sigma) \curvearrowright \tau\big)  = F_j(\sigma) \,D_a F_i(\tau);
\end{equs}
a special case of \eqref{star2_morphism}. This is a huge simplification in comparison to the original proof given in Bruned, Chandra, Chevyrev and Hairer's work \cite{BCCH}, where the authors had to go through an extended space of rooted trees in Section $4$ therein. Identity \eqref{star2_morphism} was observed in the simpler context of rough differential equation, in Lemma 3.4 of Bonnefoi, Chandra, Moinat and Weber's work \cite{BCMW}. The $\star$ product happens to be the adjoint of the Butcher-Connes-Kreimer coproduct in that setting.

\bigskip

\section{Preparation maps and multi-pre-Lie morphisms}
\label{SectionLocalProduct}

\subsection{Definition and properties}

For $\tau\in T$ denote by $\vert\tau\vert$ the number of noise symbols that appear in $\tau$. Recall from \cite{BR18} the following definition.

\smallskip

\begin{defn*}
A \textbf{\textsf{preparation map}} is a linear map $R : T\rightarrow T$ such that 
\begin{itemize}
   \item[$\bullet$] for each $ \tau \in T $ there exist finitely many $\tau_i \in T$ and constants $\lambda_i$ such that
\begin{equs} \label{analytical}
R \tau = \tau + \sum_i \lambda_i \tau_i, \quad\textrm{with}\quad \deg(\tau_i) \geq \deg(\tau) \quad\textrm{and}\quad |\tau_i| < |\tau|
\end{equs} 
   \item[$\bullet$] one has 
 \begin{equs} \label{Commutation_R}
 \left( R \otimes \id \right) \Delta = \Delta R.
 \end{equs}
\end{itemize}
\end{defn*}

\smallskip

Preparation maps are the building bricks from which renormalization maps can be constructed. A typical example of preparation map keeps fixed any tree in $\CI(T)$ and only acts non-trivially on trees with multiple edges at the root. This accounts for the fact that such trees represent products of analytical quantities, some of which needs to be renormalized to be given sense. The `deformed product' provided by $R(\tau)$ for such trees $\tau$ makes precisely that. Preparation maps were named for that reason `{\it local product renormalization maps}' in Chandra, Moinat and Weber's work \cite{CMW} for establishing a priori bounds for the $ \phi^{4}_{4-\delta}$ models in the full subcritical regime, as well as in Bruned's work \cite{Br1111} on the renormalization of branched rough paths. 

A preparation map is in particular a perturbation of the identity by elements that are more `regular' ($\deg(\tau_i) \geq \deg(\tau)$) and defined with strictly less noises ($|\tau_i| < |\tau|$). Note that the linear map $R-\textrm{Id}$ is nilpotent as a consequence of condition \eqref{analytical}. Identity \eqref{Commutation_R} encodes the fact that the recentering operator and the preparation map commute. The next statement is a direct consequence of the duality relation \eqref{def_coaction} between the product $\star$ and the splitting map $\Delta$.

\medskip

\begin{prop} \label{prop_equivalence}
Identity \eqref{Commutation_R} is equivalent to having
\begin{equs} \label{Commutation_R*}
R^{*} \left( \sigma \star \tau \right)  =  \sigma \star \left( R^{*} \tau \right)
\end{equs}
for all $\sigma \in T^+ $ and $\tau \in T$.
\end{prop}

\medskip

\begin{Dem}
We use the duality relation \eqref{def_coaction} between $\Delta$ and $\star$ to write for $\mu,\nu\in T$ and $\sigma\in T^+$
$$
\langle  \mu \otimes \sigma , \Delta R \nu  \rangle = \langle \sigma \star \mu , R \nu  \rangle = \langle R^* \left( \sigma \star \mu \right),  \nu  \rangle.
$$
The result is thus a consequence of the identity
$$
\big\langle  \mu \otimes \sigma , \left( R \otimes \textrm{Id} \right) \Delta \nu \big\rangle = \langle  R^{*} \mu \otimes \sigma ,  \Delta\nu  \rangle =  \langle  \sigma \star \left( R^*\mu \right) ,  \nu  \rangle.
$$
\end{Dem}

\medskip

\begin{defn*}
A \textbf{\textsf{strong preparation map}} is a preparation map satisfying identity \eqref{Commutation_R*} for all $\sigma\in T$ and $\tau\in T$ -- and not only for $\sigma\in T^+$. One says that $R^*$ is a right derivation for the product $\star$.
\end{defn*}

\medskip

Taking specific $\sigma$'s yields special identities. For $ \sigma = \CI_{a}(\sigma_1) $ identity \eqref{Commutation_R*} reads
\begin{equs} \label{condition_1}
R^{*}  \big( \mathcal{I}_{a}(\sigma_1) \curvearrowright \, \tau \big)  & =   \mathcal{I}_{a}(\sigma_1)\curvearrowright(  R^*\tau),
\end{equs}
that is
\begin{equs} \label{condition_11}
R^*\big(\sigma_1\curvearrowright_a \, \tau\big) = \sigma_1\curvearrowright_a( R^*\tau)
\end{equs}
Another interesting case is when $ \sigma $ is equal to the empty forest $\textbf{\textsf{1}}$ -- single node trees are identified to the empty forest when using the operator $\curvearrowright$. In that case, one has for all $k\in\N^{d+1}$
\begin{equation} \label{condition_22}
R^*(\uparrow^k\tau) = \,\uparrow^k(R^* \tau).
\end{equation}
Take care that $\uparrow^k$ in the left hand side is $\uparrow^k_{N_\tau}$, while $R^*\tau=\sum_i \tau_i$ and one has on the right hand side $\uparrow^k(R^* \tau) = \sum_i \uparrow^k_{N_{\tau_i}}(\tau_i)$. Note that the universal property of $T$ stated in Proposition \ref{freeness} implies that identities \eqref{condition_11} and \eqref{condition_22} characterize the map $R^*$ once its values on the generators $X^k\zeta_l$ are given.

\medskip

Denote by $\mcB^-$ the canonical basis of $T^-$ and recall from \cite{BHZ} or \cite{BaiHos} that $\mathbb{R}[T^-]$ is equipped with an (Hopf) algebra structure. We follow \cite{BR18} and define for any character $\ell$ of $\mathbb{R}[T^-]$ and all $\tau\in T$ 
\begin{equs} \label{def_R}
R^{*}_\ell(\tau) := \sum_{\sigma \in \mathcal{B}^-} \frac{\ell(\sigma)}{S(\sigma)} \,( \tau \star \sigma).
\end{equs}
(This definition corresponds to the dual of its usual definition -- see Corollary 4.5 in \cite{BR18}.) The BPHZ renormalization map from \cite{BHZ, CH16} corresponds to a particular choice of character $\ell$ on $T^-$.

\medskip

\begin{prop} 
The maps $ R^{*}_{\ell}$ are strong preparation maps.
\end{prop}

\medskip

\begin{Dem} 
From definition \eqref{def_R}, one has for any $\mu,\tau\in T$ 
\begin{equs}
R^*_\ell (\mu \star \tau) = \sum_{\sigma \in \mathcal{B}^-} \frac{\ell(\sigma)}{S(\sigma)} \,( \mu \star \tau) \star \sigma.
\end{equs}
By using the associativity of $\star$ one gets 
\begin{equs}
R^*_\ell( \mu \star \tau) = \sum_{\sigma \in \mathcal{B}^-} \frac{\ell(\sigma)}{S(\sigma)} \,  \mu \star (\tau\star \sigma) =  \mu \star R^{*}_\ell(\tau).
\end{equs}
The condition on the degree  $\deg(\cdot) $ in \eqref{analytical} comes from the fact that we are summing over decorated trees with negative degree in the definition of
$ R^{*}_{\ell} $. For $ |\cdot| $ which measures the size of the trees, this is a consequence of the definition of the $\star$ products, which breaks any decorated tree into two parts of smaller size. 
\end{Dem}

\medskip

It is not clear presently whether preparation maps are actually always strong. This holds however in the special case of the Butcher-Connes-Kreimer Hopf algebra, involved in the study of branched rough paths. Although elementary, the next statement will play a crucial role in the proof of our main result, Theorem \ref{ThmMain}, in the next section.

\medskip

\begin{prop} \label{PropCombined}
For every $1\leq i\leq k_0$, for every $\tau =  X^{k} \zeta_{l}\prod_{j=1}^n \CI_{a_j}(\tau_j)$ with $a_j=(\frak{t}_{l_j},p_j)\in\frak{T}^+\times\N^{d+1}$, one has
\begin{equs}
F_i(R^*\tau) = \partial^k D_{a_1} ... D_{a_n} F_i(R^*\zeta_l)\,\prod_{j=1}^{n} F_{l_j}(\tau_j)
\end{equs}
\end{prop}

\medskip

\begin{Dem}
Writing $\tau = \left( X^{k} \prod_{j=1}^{n} \CI_{a_j}(\tau_j) \right) \star \zeta_{l}$, and using the right derivatin property \eqref{Commutation_R*}, one gets
\begin{equs}
R^{*} \left( \left( X^{k} \prod_{j=1}^n \CI_{a_j}(\tau_j) \right) \star \zeta_{l} \right) = \left(  X^{k} \prod_{j=1}^n \CI_{a_j}(\tau_j) \right) \star \left( R^{*} \zeta_{l} \right) 
\end{equs}
so identity \eqref{star2_morphism} in Proposition \ref{star_morphism} yields
\begin{equs}
F_i\bigg(X^{k} \prod_{j=1}^n \CI_{a_j}(\tau_j)  \star \left( R^{*} \zeta_{l} \right) \bigg) = \partial^k   D_{a_1} ... D_{a_n} F_i(R^*\zeta_l)\,\prod_{j=1}^{n} F_{l_j}(\tau_j).
\end{equs}
\end{Dem}
 
\medskip 

\begin{defn*}
A \textbf{\textsf{good multi-pre-Lie morphism on $T$}} is a map $A : T\rightarrow T$ such that one has for all $\sigma, \tau\in T$ and $k\in\N^{d+1}$, and all $a\in\frak{L}^+\times\N^{d+1}$,
$$
A(\sigma\curvearrowright_a\tau) = (A\sigma)\curvearrowright_a(A\tau), \quad \textrm{and}\quad A\hspace{-0.05cm}\uparrow^k\,=\,\uparrow^k\hspace{-0.1cm}A.
$$
\end{defn*}

\medskip

The fundamental role of (good) multi-pre-Lie morphisms in renormalization matters was unveiled first in a rough paths setting in \cite{BCFP}, and then in \cite{BCCH}, in a regularity structures setting. Let $R$ stand for a preparation map. This is the key feature of renormalization maps that allows to obtain the renormalized equation. We associate to a strong preparation map $R$ a linear map $M^{\circ} : T\rightarrow T$, defined by the requirement that $M^{\circ} \one  = \one$, that $M^{\circ}$ is multiplicative, by the data of the $M^{\circ}\zeta_l$, for $1\leq l\leq n_0$, and the induction relation
\begin{equ}[e:defM]
M^{\circ}\big(\CI_{(\Labhom,k)}\tau\big) = \CI_{(\Labhom,k)}\big(M^\circ(R\tau)\big)
\end{equ}
for all $\tau\in T$ and $(\frak{t},k)\in\mathcal{L}^+\times\N^{d+1}$. Define, as in Section 3.1 of \cite{BR18},
\begin{equation} \label{EqConstructionRecipe}
M := M^\circ R.
\end{equation}

\medskip

\begin{thm} \label{adjoint_M} 
The map $M^*$ is a good multi-pre-Lie morphism satisfying $M^*(\zeta_l) = R^*(\zeta_l)$, for all $1\leq l\leq n_0$.
\end{thm}

\medskip

\begin{Dem} 
$\bullet$ We start by showing that the map $(M^\circ)^*$ is also multiplicative. One has
\begin{equs}
\big\langle (M^{\circ})^{*} \CI_a(\sigma),  \CI_{a}(\tau)\big \rangle &=  \big\langle  \CI_a(\sigma), M^{\circ} \CI_{a}(\tau) \big\rangle = \big\langle  \CI_a(\sigma),  \CI_{a}(M \tau) \big\rangle   \\ 
&= \langle  \sigma,  M \tau  \rangle = \langle  M^{*} \sigma,   \tau \rangle = \big \langle  \CI_a(M^{*}\sigma),  \CI_{a}(\tau) \big\rangle,
\end{equs}
which implies $(M^{\circ})^{*} \CI_a(\sigma)  =  \CI_a(M^{*}\sigma)$. Denote by $\Delta_d$ the deconcatenation coproduct
\begin{equs}
\Delta_d \CI_{a}(\tau) = \CI_{a}(\tau) \otimes \one + \one \otimes \CI_{a}(\tau), \quad \Delta_d X^k = X^k \otimes \one + \one\otimes X^k,
\end{equs}
extended multiplicatively not to the tree product but to the product between a decorated tree and a decorated tree with no polynomial decorations at the root. It follows then from the multiplicativity of $M^\circ$ and the identity $M^\circ \CI_a(\sigma)  =  \CI_a(M\sigma)$, that
\begin{equs}
\big\langle (M^{\circ})^{*} \CI_{a}(\sigma) \mu,  \tau \big\rangle = \big\langle  \CI_{a}(\sigma) \mu, M^{\circ} \tau \big\rangle = \big\langle  \CI_{a}(\sigma) \otimes \mu, \Delta_d M^{\circ} \tau \big\rangle,
\end{equs}
and
\begin{equs}
\Delta_d M^{\circ} = \left( M^{\circ} \otimes M^{\circ} \right)\Delta_d.
\end{equs}
Then, we get
\begin{equs}
\big\langle  \CI_{a}(\sigma) \otimes \mu, \Delta_d M^{\circ} \tau \big\rangle &= \big\langle  \CI_{a}(\sigma) \otimes \mu, \left( M^{\circ} \otimes  M^{\circ} \right)\Delta_d \tau \big\rangle   \\ 
&= \big\langle  (M^{\circ})^{*} \CI_{a}(\sigma) \otimes (M^{\circ})^{*} \mu, \Delta_d \tau \big\rangle = \big\langle  (M^{\circ})^{*} \CI_{a}(\sigma)  (M^{\circ})^{*} \mu,  \tau \big\rangle
\end{equs}
which concludes the proof of the multiplicativity.

\smallskip

$\bullet$ It will be useful for our purpose to decompose the grafting map $\curvearrowright_a$ into the sum of a grafting map at the root and a grafting map outside the root
$$
\curvearrowright_a = \curvearrowright_a^\textrm{root} + \curvearrowright_a^\textrm{non-root},
$$
with, for $\tau=X^k\prod_{j=1}^n\CI_{a_i}(\tau_j)$ and $a_i\in\frak{T}\times\N^{d+1}$,
$$
\sigma\curvearrowright_a^\textrm{root}\tau := \sum_{m\in\N^{d+1}} {k\choose m}\,X^{k-m}\,\CI_{a-m}(\sigma)\prod_{j=1}^n\CI_{a_j}(\tau_j)
$$
and
$$
\sigma\curvearrowright_a^\textrm{non-root}\tau := X^k \sum_{i=1}^n \CI_{a_i}(\sigma\curvearrowright_a\tau_i)\prod_{j\neq i}\CI_{a_j}(\tau_j).
$$
We proceed by induction on the size of the trees appearing in the product. In the induction hypothesis, we include the two following identities for $ \sigma, \tau \in \CT $:
 \begin{equs} \label{ident_1}
 M^{*} \left( \sigma \curvearrowright_a \tau \right) =   (M^{*} \sigma) \curvearrowright_a  (M^{*}  \tau)
 \end{equs}
 and
 \begin{equs} \label{ident_2}
 (M^{\circ})^{*} \left( \sigma \curvearrowright_a \tau \right) =   (M^{*} \sigma) \curvearrowright_a   \big((M^{\circ})^{*} \tau\big).
 \end{equs}
Let $ \sigma, \tau \in \CT $, one has
\begin{equs}
M^{*} \left( \sigma \curvearrowright_a \tau \right) & =  R^* (M^{\circ})^{*}\left( \sigma\curvearrowright_a \tau \right)   \\ 
&=  R^* \left( M^{*} \sigma \curvearrowright_a (M^{\circ})^{*} \tau \right)   \\ 
&=  \left( M^* \sigma  \curvearrowright_a R^{*}(M^{\circ})^{*} \tau \right) = \left( M^{*} \sigma\curvearrowright_a M^{*} \tau \right)
\end{equs}
where we have applied the induction hypothesis given by \eqref{ident_2} on $  (M^{\circ})^{*} $ and the righ-morphism property of $ R^{*} $. We consider $ \tau = X^k \bar \tau $ where $ \bar \tau = \prod_{i=1}^n \CI_{a_i}(\tau_i) $. The multiplicativity property of $(M^{\circ})^*$ and the fact that $(M^{\circ})^*\CI_a(\sigma)  =  \CI_a(M^*\sigma)$ yield
\begin{equs}
(M^{\circ})^{*}\left( \sigma\curvearrowright_a^{\root}\tau \right) & = (M^{\circ})^{*}  \sum_{\ell \in \N^{d+1}} {k \choose \ell} X^{k-\ell} \CI_{a-\ell}(\sigma) \bar \tau 
\\ &  =    \sum_{\ell \in \N^{d+1}} {k \choose \ell} X^{k-\ell} \CI_{a-\ell}(M^{*}\sigma) (M^{\circ})^{*} \bar \tau 
 \\ & = M^{*} \sigma\curvearrowright_a^{\root}(M^{\circ})^{*} \tau 
\end{equs}
For the grafting outside the root, we use the induction hypothesis. One has:
\begin{equs}
(M^{\circ})^{*}\left( \sigma\curvearrowright_a^{\nonroot}\tau \right)  &= (M^{\circ})^{*}  \CI_{a_i}\left( \sigma\curvearrowright_a  \tau_i \right) \prod_{j \neq i} \CI_{a_j}(\tau_j)   
\\ & =   \CI_{a_i}\left( M^{*} (\sigma\curvearrowright_a  \tau_i) \right) \prod_{j \neq i} \CI_{a_j}(M^{*} \tau_j) 
\\ & =   \CI_{a_i}\left( M^{*} \sigma\curvearrowright_a  M^{*} \tau_i \right) \prod_{j \neq i} \CI_{a_j}(M^{*} \tau_j)
\\
& = M^* \sigma \curvearrowright_a^{\nonroot}(M^{\circ})^{*}\tau.
\end{equs}
where we have used the induction hypothesis \eqref{ident_1} on $ \sigma $ and $ \tau_i $. The proof that $M^*\uparrow^k = \uparrow^{k} M^*$, works the same by decomposing insertion of polynomial decorations at the root and insertion outside the root.
\end{Dem}

\medskip

Recall the setting of branched rough paths involves the Butcher-Connes-Kreimer Hopf algebra  -- basics on branched rough paths can be found in Gubinelli's original article \cite{Gub10}, Hairer and Kelly's work \cite{HK}, or \cite{BaiSimple}, for instance.

\medskip

\begin{cor} \label{R_translation}
In the Butcher-Connes-Kreimer setting, (good) multi-pre-Lie morphisms are in bijection with preparation maps.
\end{cor}

\medskip

\begin{Dem} 
We already noticed that all preparation maps are strong in the Butcher-Conner-Kreimer setting. There is no polynomial decorations in this setting, and (automatically strong) preparation maps $R$ define (good) multi-pre-Lie morphisms $M$, with the map $R\mapsto M$ being injective because $ M^{\circ}$ is invertible. Only multi-pre-Lie morphisms make sense in the Butcher-Conner-Kreimer setting. Given a multi-pre-Lie morphism $M$, we essentially have no choice for $R$; it needs to satisfy
$$
R^*(\zeta_l) := M^*(\zeta_l), \quad R^*(\sigma\curvearrowright_a\tau)  = \sigma\curvearrowright_a( R^{*} \tau).
$$
The right-morphism property of $R^*$ gives back the `right commutation' relation \eqref{Commutation_R} with the coproduct. Property \eqref{analytical} of the map $R$ just defined can then be read off on the identity
$$
\big\langle \sigma\curvearrowright_a\tau\,,\,R\mu\big\rangle = \big\langle R^*(\sigma\curvearrowright_a\tau)\,,\,\mu\big\rangle = \sum \langle R^*\tau,\mu_1\rangle\,\langle\sigma,\mu_2\rangle,
$$
using Sweedler's notation $\Delta\mu=\sum\mu_1\otimes\mu_2$. The map $M\mapsto R$ injective, so the conclusion follows.
\end{Dem}

\medskip

\begin{rem}
Keep working in the Butcher-Connes-Kreimer setting, assuming an (automatically good) multi-pre-Lie morphism $M$ is given. Note that if one defines $M^\circ$ by the induction relations
$$
(M^\circ)^*\zeta_l = \zeta_l, \quad (M^\circ)^*(\sigma\curvearrowright_a\tau) = M^*\sigma\curvearrowright_a (M^\circ)^*\tau,
$$
then one has indeed $M=M^\circ R$. One can check by induction that $ (M^{\circ})^{*} $ is multiplicative as a consequence of the fact that $(M^\circ)^* \one = \one $ and $ (M^\circ)^* \CI_{a}(\tau) = \CI_{a}(M^* \tau) $. One uses the induction hypothesis to see that
\begin{equs}
(M^\circ)^*(\sigma\curvearrowright^{\nonroot}_a\tau) = M^*\sigma\curvearrowright^{\nonroot}_a (M^\circ)^*\tau,
\end{equs} which implies the multiplicativity given by 
 \begin{equs}
(M^\circ)^*(\sigma\curvearrowright^{\root}_a\tau) = M^*\sigma\curvearrowright^{\root}_a (M^\circ)^*\tau.
\end{equs}
So Corollary \ref{R_translation} entails that all multi-pre-Lie morphisms $M$ -- hence all renormalization maps, are obtained in that setting from a preparation map $R$ using the construction \eqref{e:defM} and \eqref{EqConstructionRecipe}. Such a statement is open for regularity structures.
\end{rem}

\bigskip

\subsection{Models associated to preparation maps}

Let kernels $(K_i)_{1\leq i\leq \ell_0}$ with a polynomial singularity at $0$, and smooth noises $(\xi_l)_{1\leq l\leq n_0}$ on the state space be given. Following \cite{BR18}, one can associate to a preparation map $R$ an admissible model $\sf M$ on $T$. It is defined from a side family $\big((\Pi_x^{M^\circ}\tau)(\cdot)\big)_{x,\tau}$ of smooth functions on the state space satisfying
$$
\big(\Pi_x^{M^{\circ}} \one\big)(y) = 1  , \quad \big(\Pi_{x}^{M^{\circ}} \zeta_l\big)(y)  = \xi_l(y) , \quad \big(\Pi_{x}^{M^{\circ}} X_i\big)(y) = y_i-x_i, 
$$
{\it the multiplicativity condition}
$$
\big(\Pi_{x}^{M^{\circ}}(\sigma\tau)\big)(y) = \big(\Pi_{x}^{M^{\circ}}\sigma\big)(y) \, \big(\Pi_{x}^{M^{\circ}}\tau\big)(y),
$$
and the condition
\begin{equation} \label{EqConditionAlmostAdmissibility}  \begin{aligned}
\Big(\Pi_{x}^{M^{\circ}} \CI_{a} (\tau)\Big)(y) = \Big( D^k K_i* \Pi^{M^\circ}_{x}(R\tau) \Big)(y) - \sum_{|\ell|_{\s} \leq \deg( \CI_{a} (\tau) )} \frac{(y-x)^{\ell}}{\ell!} \Big(D^{k + \ell} K_i* \Pi_{x}^{M^\circ}(R\tau)\Big)(x),
\end{aligned} \end{equation}
for $a=(\frak{t}_i,k)$. Define for all $x$ and $\tau$ a smooth function on the state space
$$
\big(\Pi^{(R)}_x\tau\big)(\cdot) := \Big(\Pi^{M^\circ}_x(R\tau)\Big)(\cdot).
$$
Bruned gave in Proposition 3.16 of \cite{BR18} an explicit contruction of an admissible model $\sf M=(g,\Pi)$ on $T$, with values in the space of smooth functions, such that the operators $\Pi^{(R)}_x$ are indeed associated with $\sf M$, in the sense that one has for all $\tau\in T$ and $x$
$$
\Pi^{(R)}_x\tau = {\sf \Pi}_x^{\sf g}\tau.
$$
Since the model $\sf M$ takes values in the space of continuous functions, the reconstruction operator $\sf R^M$ associated with it is given by the explicit formula
$$
\big({\sf R^M v}\big)(x) = \big({\sf \Pi}_x^{\sf g}{\sf v}(x)\big)(x)
$$
for any modelled distribution $\sf v$ with positive regularity, so 
$$
\big({\sf R^M v}\big)(x) = \Big(\Pi_x^{M^\circ}R{\sf v}(x)\Big)(x).
$$
Emphasize that $M^\circ$ depends on $R$, so the operators $\Pi_x^{M^\circ}$ are different from the naive interpretation operators one obtains when $R=\textrm{Id}$. For preparation maps $R$ for which 
\begin{equation} \label{EqConditionR}
R(\CI_a\tau) = \CI_a\tau
\end{equation}
for all $\tau\in T$ and $a\in\frak{L}^+\times\N^{d+1}$, one has 
$$
{\sf \Pi}_x^{\sf g}(\CI_a\tau) = \Pi_x^{M^\circ}(\CI_a\tau).
$$
(Condition \eqref{EqConditionR} is consistent with the idea that the preparation/`local product' map $R$ `renormalizes' only ill-defined products -- no product is involved at the root of the tree $\CI_a\tau$.) The point here is that $\Pi_x^{M^\circ}$ is multiplicative while ${\sf \Pi}_x^{\sf g}$ is not. Denote by $T_X\subset T$ the linear space spanned by polynomial in $T$. Modelled distribution $\sf v$ with values in the subspace $\CI(T)\oplus T_X$ of $T$ satisfy in that case the identity
\begin{equation} \label{EqRelationPiMPiMCirc}
\big({\sf R^M v}\big)(x) = \Big(\Pi_x^{M^\circ}{\sf v}(x)\Big)(x).
\end{equation}
It follows further from relation \eqref{EqConditionAlmostAdmissibility} that the model $\sf M$ is admissible if condition \eqref{EqConditionR} holds.

\bigskip

\section{A short proof for the renormalised equation}
\label{SectionShortProof}

This section contains the statement and proof of our main result, Theorem \ref{ThmMain}, describing the autonomous dynamics satisfied by ${\sf R^Mu}$ when $\sf M$ is the model constructed from a preparation map $R$ and $\sf u$ is the solution to the lift of system \eqref{main_equation} to its associated regularity structure 
\begin{equation} \label{EqLiftedSystem}
{\sf u}_i = \mathcal{K}^{\sf M}_i\Big(\mathcal{Q}_{\gamma-2} \big({\sf F}_i({\sf u}, D{\sf u})\zeta\big)\Big) + \mathcal{P}_\gamma u_i(0), \qquad (1\leq i\leq k_0).
\end{equation}
The operator $\mathcal{P}_\gamma$ stands here for the projection on the subspace of elements of degree less than $\gamma$ of the natural lift in the polynomial regularity structure $T_X$ of the map $x=(t,x')\mapsto \big(P_tu(0)\big)(x')$.  (Recall $x$ stands for a generic spacetime point.) The notation $\mathcal{Q}_{\gamma-2}$ stands here for the natural projection from $T$ to the linear subspace $T_{<\gamma-2}$ of elements of $T$ of degree less than $\gamma-2$. The $\sf M$-dependent map $\mathcal{K}^{\sf M}_i$ is the regularity structure lift of the operator $K_i=(\partial_t-L_i)^{-1}$; it sends continuously the space $\mathcal{D}^{\gamma-2, \eta}(T,\sf g)$ into $\mathcal{D}^{\gamma, \eta'}(T,\sf g)$. (The exponent $\beta$ in the spaces $\mathcal{D}^{\alpha, \beta}(T,\sf g)$ is related to the behaviour of the function near time $0^+$. We refer the reader to \cite{BHZ,BCCH} or \cite{BaiHos} for a full account.) From the definition of the operators $\mathcal{K}^{\sf M}_i$, for a solution ${\sf u}=({\sf u}_1,\dots,{\sf u}_{k_0})$ of system \eqref{EqLiftedSystem}, each map ${\sf u}_i$ takes values in $\CI_{(\frak{t}_i,0)}(T)\oplus T_X$. Write
$$
{\sf u}_i =: \sum {\sf u}_{i,\tau}\tau,
$$
for a sum over trees $\tau$ in the canonical basis of $\CI_{(\frak{t}_i,0)}(T)\oplus T_X$ -- monomials are seen as trees with just one vertex here. We recall here from Section 4.2 of \cite{reg} the definition of the lift ${\sf F}_i$ of the smooth enough function $F_i$. One has for any ${\sf a} =: a_{\bf 1}{\bf 1} + {\sf a}'\in T$, with $\langle{\sf a}',{\bf 1}\rangle=0$,
\begin{equation} \label{EqDefnRSLiftFunction}
{\sf F}_i({\sf a}) = \sum_k\frac{D^kF(a_{\bf 1})}{k!}\,({\sf a}')^k.
\end{equation}
One of the main results of \cite{BCCH} states that for all $\tau\in T$ in the canonical basis, with degree less than $\gamma-2$, 

\begin{equation} \label{EqCoherence}
{\sf u}_{i,\tau} =  \frac{{\sf F}_i(\tau)({\sf u}, D{\sf u})}{S(\tau)}.
\end{equation}
(Note that ${\sf F}_i(\tau)=0$ for all trees $\tau\in T$ that are not generated by the rule associated with the non-linearities $F_i$ in \eqref{main_equation}.) Let $\xi=(\xi_1,\dots,\xi_{n_0})$ stand for a smooth function.

\medskip

\begin{thm} \label{ThmMain}
Let $R$ be a strong preparation map such that 
$$
R\tau=\tau, \qquad\textrm{for} \quad \tau\in\big(\CI(T)\oplus T_X\big).
$$ 
Let $\sf M$ stand for its associated admissible model. Then $u:={\sf R^Mu}$ is a solution of the renormalized system
\begin{equation} \label{EqRenormalizedSystem}
(\partial_t - L_i) u_i = F_i(u,\nabla u)\,\xi + \sum_{l=1}^{n_0} F_i\Big(\big(R^* - \textrm{\emph{Id}}\big) \zeta_l\Big)(u,\nabla u)\,\xi_l, \qquad (1\leq i\leq k_0).
\end{equation}
\end{thm}

\medskip

\begin{Dem}
As we are working with an admissible model we have
$$
(\partial_t-L_i)u_i = {\sf R^M}\Big(\mathcal{Q}_{\gamma-2}\big({\sf F}_i({\sf u}, D{\sf u})\zeta\big)\Big) =: {\sf R^M}({\sf v}_i)
$$
with
$$
{\sf v}_i = \sum_{\textrm{deg}(\tau)<\gamma-2} \frac{{\sf F}_i(\tau)({\sf u}, D{\sf u})}{S(\tau)} \,\tau,
$$
for a sum over the canonical basis of $T$, from the coherence condition \eqref{EqCoherence}. The function ${\sf v}_i$ is a modelled distribution of regularity $\gamma$. One has by construction
$$
\big({\sf R^M}{\sf v}_i\big)(x) = \Big(\Pi_x^{M^\circ}R{\sf v}_i(x)\Big)(x),\qquad \textrm{and}\qquad \langle R{\sf v}_i,\tau\rangle = \langle {\sf v}_i\,,\,R^*\tau\rangle = F_i(R^*\tau),
$$
with
$$
F_i(R^*\tau) = \partial^kD_{a_1}\cdots D_{a_n} F_i(R^*\zeta_l) \prod_{j=1}^nF_{l_j}(\tau_j)
$$
when $\tau=X^k\zeta_l\prod_{j=1}^n\CI_{a_j}(\tau_j)$ and $a_j=(\frak{t}_{n_j},k_j)$, from Proposition \ref{star_morphism}. One can thus rewrite the equality
$$
R{\sf v}_i = \sum \frac{F_i(R^*\tau)}{S(\tau)}\,\tau
$$ 
under the form
\begin{equation*}
R{\sf v}_i = \sum_{l,n} \sum_{a_1,...,a_n} \sum_{k}\frac{k! \prod_{i=1}^n S(\tau_i)}{S\big(X^k\zeta_l\prod_{j=1}^n\CI_{a_j}(\tau_j)\big)} \, \frac{X^k}{k!} \prod_{j=1}^n \sum_{\tau_j\in \CT} \frac{F_{l_j}(\tau_j)}{S(\tau_j)}\, \CI_{a_j}(\tau_i) \, \partial^k \prod_{i=1}^n D_{a_i}F_i(R^*\zeta_l) \zeta_l.
\end{equation*}
chopping $\tau=X^k\zeta_l\prod_{j=1}^n\CI_{a_j}(\tau_j)$ in different pieces. Using distinct $a_j$'s
\begin{equs}
R{\sf v}_i= \sum_{l,n} \sum_{a_1,...,a_n} \sum_{k} \frac{X^k}{k!} \prod_{j=1}^n \left( \sum_{\tau_j \in T} \frac{1}{\beta_j !} \frac{F_{l_j}(\tau_j)}{S(\tau_j)}\,\CI_{a_j}(\tau_j) \right)^{\beta_j} \partial^k \prod_{j'=1}^n (D_{a_{j'}})^{\beta_{j'}} F_{{l_j{'}}} (R^*\zeta_l) \zeta_l,
\end{equs}
and the Fa\`a di Bruno formula from Lemma A.1 in \cite{BCCH}
\begin{equs}
\frac{\partial^{k}G}{k!} =  \sum_{b_1,...,b_m} \sum_{k = \sum_{j=1}^m \beta_j k_j} \prod_{j=1}^m \frac{1}{\beta_j !} \left(\frac{Z_{b_i + k_i}}{k_j!}\right)^{\beta_j} \prod_{j=1}^m (D_{b_j})^{\beta_j} G 
\end{equs}
for a function $G$ of distinct variables $Z_{b_1},\dots, Z_{b_m}$, one obtains
\begin{equs}
R{\sf v}_i = \sum_{l,n} \sum_{a_1,...,a_n} \sum_{\beta_1,...,\beta_n } \prod_{j=1}^n \frac{1}{\beta_i !} \left({\sf u}_{a_i} - \langle {\sf u}_{a_j}, \one \rangle \right)^{\beta_j} \prod_{j'=1}^n (D_{a_{j'}})^{\beta_{j'}} F_i (R^*\zeta_l) \zeta_l.
\end{equs}
From the definition of ${\sf F}_i$ recalled in \eqref{EqDefnRSLiftFunction}, this is equivalent to

\begin{equs} \label{main_identity}
R{\sf v}_i = \sum_{l=1}^{n_0} {\sf F}_i(R^{*}\zeta_l)({\sf u}, D{\sf u}) \zeta_l.
\end{equs}
Recall $x$ stands for a generic spacetime point. Using the (crucial) multiplicativity of $\Pi_x^{M^\circ}$ and identity \eqref{EqRelationPiMPiMCirc} giving back $\big({\sf R^Mu}\big)(x)$ in terms of $\Pi_x^{M^\circ}$, we see that 
$$
\big((\partial_t-L_i)u_i\big)(x) = \big({\sf R^M}{\sf v}_i\big)(x)  = \Pi_x^{M^\circ}\big(R{\sf v}_i(x)\big)(x)
$$

\begin{equation*} \begin{split}
&= \sum_{l=1}^{n_0} \Pi_x^{M^\circ}\Big({\sf F}_i\big(R^{*}\zeta_l)({\sf u}(x), D{\sf u}(x)\big) \zeta_l\Big)(x)   \\
&= \sum_{l=1}^{n_0} F_i(R^*\zeta_l)\Big(\big(\Pi_x^{M^\circ}{\sf u}(x)\big)(x), \nabla \big(\Pi_x^{M^\circ}{\sf u}(x)\big)(x)\Big) \,\Pi_x^{M^\circ}\zeta_l   \\
&= \sum_{l=1}^{n_0} F_i(R^*\zeta_l)\big(u(x),\nabla u(x)\big)\xi_l.
\end{split}\end{equation*} 
\end{Dem}

\medskip

In the particular case where the preparation map $R$ is of BPHZ form \eqref{def_R}, a direct computation shows that the renormalized system \eqref{EqRenormalizedSystem} takes the form \eqref{EqRenormalisedEquation} if we further assume that $R_{\ell^{\varepsilon}}^*\zeta_l=\zeta_l$ for all $l\neq 0$ -- recall $\zeta_0=\textbf{\textsf{1}}$, since one has from  \eqref{def_R}
\begin{equs}
F_i(R_{\ell^{\varepsilon}}^* \one  -\one) = \sum_{\tau \in \mathcal{B}^{-}\setminus \{\one \} } \ell^{\varepsilon}(\tau)\frac{F_i(\tau)}{S(\tau)}.
\end{equs}
This assumption accounts for the fact that we never need to substract a multiple of one of the noises $(\zeta_l)_{1\leq l\leq n_0}$ to any tree-indexed quantity in our renormalization algorithm. Note that Chandra, Moinat and Weber used in \cite{CMW} a similar strategy to get back the renormalized equation in the particular case of the $\Phi^4_{4-\delta}$ equation.

\bigskip
\bigskip

\noindent \textcolor{gray}{$\bullet$} {\sf I. Bailleul} -- Univ. Rennes, CNRS, IRMAR - UMR 6625, F-35000 Rennes, France.   \\
\noindent {\it E-mail}: ismael.bailleul@univ-rennes1.fr   

\medskip

\noindent \textcolor{gray}{$\bullet$} {\sf Y. Bruned} --  School of Mathematics, University of Edinburgh, EH9 3FD, Scotland   \\
{\it E-mail}: Yvain.Bruned@ed.ac.uk

\end{document}